\documentstyle[preprint,aps]{revtex}

\newcommand{\nc}{\newcommand}
\nc{\al}{\alpha}
\nc{\ga}{\gamma}
\nc{\de}{\delta}
\nc{\ep}{\epsilon}
\nc{\ze}{\zeta}
\nc{\et}{\eta}
\nc{\ka}{\kappa}
\nc{\la}{\lambda}
\nc{\rh}{\rho}
\nc{\si}{\sigma}
\nc{\ta}{\tau}
\nc{\up}{\upsilon}
\nc{\ph}{\phi}
\nc{\ch}{\chi}
\nc{\ps}{\psi}
\nc{\om}{\omega}
\nc{\Ga}{\Gamma}
\nc{\De}{\Delta}
\nc{\La}{\Lambda}
\nc{\Si}{\Sigma}
\nc{\Up}{\Upsilon}
\nc{\Ph}{\Phi}
\nc{\Ps}{\Psi}
\nc{\Om}{\Omega}
\nc{\ptl}{\partial}
\nc{\del}{\nabla}
\nc{\be}{\begin{eqnarray}}
\nc{\ee}{\end{eqnarray}}
\nc{\cj}{{\cal J}}
\nc{\ck}{{\cal K}}
\nc{\ov}{\overline}

\begin{document}
\title{A $q$-Lorentz Algebra From $q$-Deformed Harmonic Oscillators}
\draft
\preprint{UM-P-95/83; RCHEP-95/20}
\author{A. Ritz\footnote{ritz@physics.unimelb.edu.au}
and G. C. Joshi\footnote{joshi@physics.unimelb.edu.au}}
\address{Research Centre for High Energy Physics, \\
	 School of Physics, The University of Melbourne, \\
	 Parkville, VIC 3052 Australia}
\date{4$^{th}$ September 1995}

\maketitle

\begin{abstract}
A mapping between the operators of the bosonic
oscillator and the Lorentz rotation and boost generators is presented.
The analog of this map in the $q$-deformed regime is then applied to
$q$-deformed bosonic oscillators to generate a $q$-deformed Lorentz algebra,
via an inverse of the standard chiral decomposition.
A fundamental representation, and the co-algebra structure,
are given, and the
generators are reformulated into $q$-deformed rotations and boosts.
Finally, a relation between the $q$-boson operators and a basis of $q$-deformed
Minkowski coordinates is noted.
\end{abstract}

\newpage
\section{Introduction}
The existence of symmetries has proved to be a universal feature in many
aspects of physics. Lie group theory, as a
mathematical description of these symmetries,
has also proved to be a powerful tool for prediction and unification in the
last
thirty years.
This may suggest that consideration of possible
deformations of these symmetries would not be highly profitable. However,
it is important to realise that arbitrary deformations can be considered
as generalisations to the standard results, and can
lead to new and more powerful
symmetries that were not previously apparent.
Indeed deformations of physical theories have been
very important in the advance of
physics in the last hundred years.
For example, special relativity can be considered as a
deformation of Galilean relativity with deformation parameter $1/c$,
while quantum mechanics can similarly be regarded as a deformation of
classical mechanics with parameter $\hbar$ \cite{carow91}.

There has recently been considerable interest in the deformation, via an
arbitrary parameter $q$, of the Lie group structures commonly found in
physics. These {\em quantum groups}
\cite{drinfeld86,jimbo85,jimbo86a,jimbo86b,woron87}
are deformations of the universal enveloping algebra of the underlying
Lie group and have a Hopf algebra structure.
While quantum groups associated with many of the simple Lie algebras
initially received most attention,
deformations of the spacetime symmetry groups, i.e. the Lorentz group,
have also been considered recently by a number of authors
\cite{carow90,schmidke91,ogiev91,carow91,ogiev92,song92b,aneva94,dabrow94,pillin94}.

Ogievetsky et al. \cite{ogiev91}
have considered a chiral decomposition of a $q$-Lorentz group into
two inequivalent $SL_q(2)$ groups acting on dotted and undotted spinors.
This resulted in a {\em $q$-decomposition\/} since the two $SL_q(2)$ groups
only $q$-commuted. In
this letter we consider the construction of a $q$-deformed Lorentz
algebra which admits a complete decomposition into two commuting chiral
sub-algebras. We base this construction on a $q$-deformation
of the classical group isomorphism $SO(4) \cong SU(2) \otimes SU(2)$.
Recognising that the non-compact Lorentz group $SO(3,1)$ is equivalent to
its compact relation $SO(4)$, with an
altered metric \cite{basu77}, we can
thus consider the $q$-deformed Lorentz algebra resulting from the direct
product of two copies of the $SU_q(2)$ (or $U_q(su(2))$) quantum group in
the Drinfel'd-Jimbo basis. This algebra can then be linked back to the
$q$-boson ($q$-SHO) operators by
making use of the $q$-deformed Jordan-Schwinger map
\cite{biedenharn89,macfar89,sun89,kulish90}.

In Section 2 the mapping between harmonic oscillator annihilators and
creators and the ``classical'' Lorentz
rotations and boosts is considered. An analogous construction
is then applied to the $q$-oscillator operators to generate a $q$-deformed
Lorentz algebra in Section 3. Fundamental representations are obtained
and the generators are reformulated into $q$-deformed rotations and boosts.
It is also noted that there exists a mapping from $q$-boson operators
to the $q$-Minkowski coordinates of Ogievetsky et al. \cite{ogiev91,ogiev92}.

\section{Construction of the Lorentz group generators}
Consider two commuting sets of SHO annihilators and creators
\be
 [a_i, a_j^{\dagger}] & = & \de_{ij}, \;\;\;\;\;\; \mbox{for }i,j\in 1,2,
\ee
with number operator
\be
 N_i & = & a_i^{\dagger}a_i, \;\;\;\;\;\;\;\; i=1,2.
\ee
Then we construct three generators
\be
 J_{+} & = & a_1^{\dagger}a_2 \\
 J_{-} & = & a_2^{\dagger}a_1 \\
 J_3 & = & \frac{1}{2}(N_1 - N_2)
\ee
according to the Jordan-Schwinger map. These generators then
satisfy the commutation relations
\be
 \left[J_3, J_{\pm}\right] & = & J_{\pm} \\
 \left[J_{+},J_{-}\right] & = & 2J_3,
\ee
and generate $SU(2)$. Writing $J_{\pm}=J_1\pm iJ_2$ we have
\be
 \left[J_i,J_j\right] & = & i\ep_{ijk}J_k,
\ee
where
\be
 J_1 = \frac{1}{2}(a_1^{\dagger}a_2 + a_2^{\dagger}a_1) \\
 J_2 = \frac{1}{2i}(a_1^{\dagger}a_2 - a_2^{\dagger}a_1) \\
 J_3 = \frac{1}{2}(a_1^{\dagger}a_1 - a_2^{\dagger}a_2)
\ee
Choosing two commuting sets of $SU(2)$ generators
\be
 \left[J^m_i,J^p_j\right] & = & i \de^{mp}\ep_{ijk} J_k^p,
\ee
where $p,m=1,2$, and $i,j,k=1,2,3$, we construct the generators
\be
 \cj_i & = & J_i^1 + J_i^2 \;\;\;\;\;\;\;\;\;\;\;\;\;\;\;\mbox{(rotations)} \\
 \ck_i & = & -i(J_i^1-J_i^2)\;\;\;\;\;\;\;\;\;\;\mbox{(boosts)}
\ee
which satisfy
\be
 \left[\ck_i,\ck_j\right] & = & -i \ep_{ijk}\cj_k \label{l1}\\
 \left[\cj_i,\cj_j\right] & = & i \ep_{ijk}\cj_k \\
 \left[\cj_i,\ck_j\right] & = & i \ep_{ijk}\ck_k.\\
\ee
These clearly generate
\be
 SO(3,1) & \cong & SU(2) \otimes SU(2),
\ee
which with the Minkowski metric becomes the Lorentz group \footnote{A change of
sign in
the right hand side of Eq. \ref{l1} will give the generator algebra of
the compact group $SO(4)$.}.
Thus combining
the two maps we can write the Lorentz group generators in terms of the
bosonic oscillator operators as follows,
\be
 \cj_1 & = & \frac{1}{2} \sum_{i=1}^2
		     (a_1^{i\dagger}a_2^i+a_2^{i\dagger}a_1^i) \label{r1}\\
 \cj_2 & = & \frac{1}{2i} \sum_{i=1}^2
		     (a_1^{i\dagger}a_2^i-a_2^{i\dagger}a_1^i) \label{r2}\\
 \cj_3 & = & \frac{1}{2} \sum_{i=1}^2
		     (N_1^i-N_2^i) \label{r3}\\
 \ck_1 & = & -\frac{i}{2} (a_1^{1\dagger}a_2^1+a_2^{1\dagger}a_1^1
		     -a_1^{2\dagger}a_2^2-a_2^{2\dagger}a_1^2) \label{b1}\\
 \ck_2 & = & -\frac{1}{2} (a_1^{1\dagger}a_2^1-a_2^{1\dagger}a_1^1
		     -a_1^{2\dagger}a_2^2+a_2^{2\dagger}a_1^2) \label{b2}\\
 \ck_3 & = & -\frac{i}{2} (N_1^1-N_2^1
		     -N_1^2+N_2^2), \label{b3}
\ee
where the lower index corresponds to the individual SHO operators in each
of the commuting copies of $SU(2)$ denoted by the upper index.

The rotation and boost operators are often
combined into one tensor $\si_{\mu\nu}$, as follows,
\be
 (\si_{\mu\nu}) & = & 2\left(\begin{array}{cc}
	  0          &    \vec{\ck}^T  \\
	  -\vec{\ck} &    (\cj_{ij})
		    \end{array} \right),
\ee
where $\vec{\ck}$ is the boost vector,
\be
 \vec{\ck}^T & = & \left( \ck_1,\ck_2, \ck_3 \right),
\ee
and $(\cj_{ij})$ is the standard $3 \times 3$ rotation matrix where
\be
 \cj_i & = & \frac{1}{2}\ep_{ijk}\cj_{jk}.
\ee
$\si^{\mu\nu}$ gives the reducible spinorial representation
of $SO(3,1)$ (or $SO(4)$) which is most naturally formulated in terms of the
elements of a Dirac algebra as follows,
\be
 & & \si^{\mu\nu}  =  \frac{i}{2} [\ga^{\mu},\ga^{\nu}],
		 \;\;\;\;\;\;\;\mbox{where} \label{spinorep}\\
 & & \{\ga^{\mu},\ga^{\nu}\}  =  2\et^{\mu\nu}.
\ee
If $\et^{\mu\nu}$ is replaced by the Euclidean metric $\de^{\mu\nu}$,
as in the case
of $SO(4)$, the Dirac algebra becomes the standard Clifford algebra.
Since the generators $\si^{\mu\nu}$
can be related to the bosonic oscillator operators this suggests that
when considered as operators
we can represent
elements of the $\ga$-matrices in terms of the bosonic annihilators
and creators. i.e. in the Weyl representation,
in $2 \times 2$ block form,
\be
 (\ga^1)^{01} & = & -(\ga^1)^{10} =
	  (a_1^{\dagger}a_2 + a_2^{\dagger}a_1)   \\
 (\ga^2)^{01} & = & -(\ga^2)^{10} =
	  -i(a_1^{\dagger}a_2 - a_2^{\dagger}a_1)   \\
 (\ga^3)^{01} & = & -(\ga^3)^{10} =
	  -(a_1^{\dagger}a_1 - a_2^{\dagger}a_2)
\ee

Note that the standard form for $\ga^5$ still holds and Dirac bi-spinors $\ps$
transforming under the spinorial representation of the Lorentz group, i.e.
$\ps \rightarrow \ps'= \exp(-i/4\theta_{\mu\nu}\si^{\mu\nu})$, can still
be decomposed into $SU(2)$ transforming spinors
\footnote{The standard chiral decomposition
of Lorentz bi-spinors produces spinors transforming
under the self representation and complex conjugate self representation of
$SL(2,{\boldmath C})$, with undotted and dotted indices respectively.
However, since we are working for convenience with
the compact $SO(4)$ group, a discussion involving $SU(2)$ spinors
is more relevant.}
via the standard chiral
projectors, $P_{\pm}=(1\pm \ga^5)/2$. Vectors and higher spin representations
obtained via direct products of the chiral spinors, $\ps_R = (1/2,0)=P_+ \ps$,
and $\ps_L=(0,1/2)=P_- \ps$, will also transform in the standard manner and
thus the formulation in terms of the bosonic oscillator operators is
entirely equivalent to the standard one.

\section{Construction of the {$q$}-Lorentz algebra}
Having introduced a mapping between the bosonic oscillator operators and the
Lorentz group generators we can consider its $q$-deformation, and thereby
obtain a q-deformation of the Lorentz algebra. Our mapping will now be
\be
    q\mbox{-SHO} & \;\; \longrightarrow  \;\;& SU_q(2) \;\;\;
                   \longrightarrow \;\;\; SO_q(3,1).
\ee
The standard $q$-oscillator is defined via the deformed commutation relations
\cite{macfar89},
\be
 \left[a,a^{\dagger}\right]_{q^{-1}} & = & aa^{\dagger}-q^{-1}a^{\dagger}a
	      = q^{N_q} \label{js1}\\
 \left[a,N_q\right] & = & a \\
 \left[a^{\dagger},N_q\right] & = & -a^{\dagger},
\ee
where $N_q$ is the $q$-number operator. The substitution $q \rightarrow q^{-1}$
also implies
\be
 aa^{\dagger}-qa^{\dagger}a
	      = q^{-N_q}. \label{js2}
\ee
The Jordan - Schwinger map can be applied as before, giving the so-called
$q$-Jordan-Schwinger map (see \cite{biedenharn89,macfar89,sun89,kulish90}).
Taking a pair of commuting $q$-SHO, $a_i,a_i^{\dagger},N_{qi}$, $i$=1,2,
we define
\be
 J_{+} & = & a_1^{\dagger}a_2 \\
 J_{-} & = & a_2^{\dagger}a_1 \\
 J_3 & = & \frac{1}{2} \left( N_{q1}-N_{q2}\right).
\ee
These are the generators of the quantum universal enveloping algebra of
$SU(2)$,
$U_q(su(2))$. i.e.
\be
 \left[J_3,J_{\pm}\right] & = & \pm J_{\pm} \label{uqsu21}\\
 \left[J_{+},J_{-}\right] & = & \frac{q^{2J_3}-q^{-2J_3}}{q-q^{-1}}
	= \left[2J_3\right]_q \label{uqsu22},
\ee
where Eq. \ref{uqsu21} follows directly, and Eq. \ref{uqsu22} follows from
consideration of Eq.s \ref{js1} and \ref{js2}.

We shall return to this formulation later in the consideration of $q$-rotation
and $q$-boost generators. However, for construction of the $q$-Lorentz algebra
in analogy to the classical chiral decomposition, a different basis for
$SU_q(2)$ is more appropriate.
If we choose a basis consisting of generators
\be
 J_+ & = & a_1^{\dagger}a_2 \\
 J_- & = & a_2^{\dagger}a_1 \\
 q^{J_3} & = & q^{1/2(N_{q1}-N_{q2})},
\ee
then the $SU_q(2)$ algebra has the form (see Appendix A):
\be
 q^{J_3}J_{\pm}q^{-J_3} & = & q^{\pm 1} J_{\pm} \\
 \left[J_+,J_-\right] & = & [2J_3]_q.
\ee
In analogy to the classical direct product, $SO(3,1) \cong SU(2)
\otimes SU(2)$, we can
construct the $q$-Lorentz algebra from two copies of $SU_q(2)$. The basis
will thus be $\{J_{\pm},q^{J_3},\ov{J}_{\pm},q^{\ov{J}_3}\}$, where
\be
 q^{J_3}J_{\pm}q^{-J_3} & = & q^{\pm 1} J_{\pm} \\
 \left[J_+,J_-\right] & = & [2J_3]_q \\
 q^{\ov{J}_3}\ov{J}_{\pm}q^{-\ov{J}_3} & = & q^{\pm 1} \ov{J}_{\pm} \\
 \left[\ov{J}_+,\ov{J}_-\right] & = & [2\ov{J}_3]_q,
\ee
where barred and unbarred generators commute. The explicit commutation
properties of the two copies of $SU_q(2)$ are in contradistinction to the
chiral decomposition considered by Ogievetsky et al. \cite{ogiev91,ogiev92},
where the two chiral subalgebras only explicitly $q$-commute.

Since the two algebras, $SU_q(2)$ and $\ov{SU}_q(2)$,
are distinct, and in line with the direct product interpretation, we note
that the fundamental representation must be four dimensional. This can be
explicitly obtained by embedding the $SU_q(2)$ fundamental representation
(see Appendix A) into $4 \times 4$ matrices as follows,
\be
 \rh(q^{J_3}) & = & \left(\begin{array}{cc}
     \left(\begin{array}{cc}
	q^{1/2}       &      0    \\
	0             & q^{-1/2}
	   \end{array}\right) &        0_2     \\
     0_2                      &        1_2
			  \end{array}\right) \\
 \rh(J_+) & = & \left(\begin{array}{cc}
     \left(\begin{array}{cc}
	0       &      1    \\
	0       &      0
	   \end{array}\right) &        0_2     \\
     0_2                      &        0_2
			  \end{array}\right) \\
 \rh(J_-) & = & \left(\begin{array}{cc}
     \left(\begin{array}{cc}
	0       &      0    \\
	1       &      0
	   \end{array}\right) &        0_2     \\
     0_2                      &        0_2
			  \end{array}\right) \\
 \rh(q^{\ov{J}_3}) & = & \left(\begin{array}{cc}
     1_2                      &        0_2     \\
     0_2                      & \left(\begin{array}{cc}
	                         q^{1/2}       &      0    \\
	                         0             & q^{-1/2}
	                               \end{array}\right)
			       \end{array}\right) \\
 \rh(\ov{J}_+) & = & \left(\begin{array}{cc}
     0_2                      &        0_2     \\
     0_2                      & \left(\begin{array}{cc}
	                         0       &      1    \\
	                         0       &      0
	                               \end{array}\right)
			       \end{array}\right) \\
 \rh(\ov{J}_-) & = & \left(\begin{array}{cc}
     0_2                      &        0_2     \\
     0_2                      & \left(\begin{array}{cc}
	                         0       &      0    \\
	                         1       &      0
	                               \end{array}\right)
			       \end{array}\right).
\ee
The barred and unbarred algebras
then form distinct $2 \times 2$ blocks, as expected for the direct
product decomposition. We also note that a similar representation has been
obtained by Dabrowski et al. \cite{dabrow94}.

The coalgebra structure of this $q$-Lorentz algebra can be obtained by
analogy with the results of Dabrowski et al. \cite{dabrow94}. The
co-product ($\De$), the co-unit ($\ep$), and the antipode ($S$) are given by
\be
 \De(q^{J_3}) & = & q^{J_3} \otimes q^{J_3} \\
 \De(J_+) & = & J_+ \otimes q^{-J_3}q^{\ov{J}_3}
			    + q^{J_3}q^{-\ov{J}_3} \otimes J_+ \\
 \De(J_-) & = & J_- \otimes q^{-J_3}q^{-\ov{J}_3}
			    + q^{J_3}q^{\ov{J}_3} \otimes J_- \\
 \De(q^{\ov{J}_3}) & = & q^{\ov{J}_3} \otimes q^{\ov{J}_3} \\
 \De(\ov{J}_+) & = & \ov{J}_+ \otimes q^{-J_3}q^{\ov{J}_3}
			    + q^{J_3}q^{-\ov{J}_3} \otimes \ov{J}_+ \\
 \De(\ov{J}_-) & = & \ov{J}_- \otimes q^{J_3}q^{\ov{J}_3}
			    + q^{-J_3}q^{-\ov{J}_3} \otimes \ov{J}_- \\
 \mbox{} \nonumber\\
 \ep(q^{J_3}) & = & \ep(q^{\ov{J}_3}) = 1 \\
 \ep(J_{\pm}) & = & \ep(\ov{J}_{\pm}) = 1 \\
 \mbox{} \nonumber\\
 \mbox{} \nonumber\\
 S(q^{J_3}) & = & q^{-J_3} \\
 S(J_{\pm}) & = & -q^{\mp 1}J_{\pm} \\
 S(q^{\ov{J}_3}) & = & q^{-\ov{J}_3} \\
 S(\ov{J}_{\pm}) & = & -q^{\pm 1}\ov{J}_{\pm}.
\ee
This completes the Hopf algebra structure of the quantum Lorentz group in this
formulation.

To obtain a formalism of the $q$-Lorentz algebra corresponding to the
standard formulation in terms of rotations and boosts we return to the
original Drinfel'd-Jimbo basis of $U_q(su(2))$.
Again choosing two commuting sets
of $U_q(su(2))$ generators, $\{J_3,J_{\pm}\}$, and
$\{\ov{J}_3,\ov{J}_{\pm}\}$, we construct the generators
\be
 \cj_1 & = & \frac{1}{2}\left(J_+ +J_- + \ov{J}_+ + \ov{J}_-\right) \\
 \cj_2 & = & \frac{1}{2}\left(J_+ - J_- + \ov{J}_+ - \ov{J}_-\right) \\
 \cj_3 & = & J_3 + \ov{J}_3 \\
 \ck_1 & = & \frac{1}{2i}\left(J_+ +J_- - \ov{J}_+ - \ov{J}_-\right) \\
 \ck_2 & = & -\frac{1}{2}\left(J_+ - J_- - \ov{J}_+ + \ov{J}_-\right) \\
 \ck_3 & = & -i\left(J_3 - \ov{J}_3\right),
\ee
where we associate $\cj_i,\; i=1,2,3$ with $q$-rotations, and
$\ck_i,\; i=1,2,3$ with $q$-boosts. With this construction these $q$
generators have exactly the same structure in terms of $q$-boson operators,
as do the classical rotations and boosts in terms of the boson operators. i.e.
Eq.s \ref{r1}-\ref{b3}
still hold in the $q$-deformed regime.

In order to make contact with the analysis of Schmidke et al. \cite{schmidke91}
and Ogievetsky et al. \cite{ogiev91,ogiev92} we can
transform our formulation of the
$q$-Lorentz generators to correspond to the Woronowicz basis of $SU_q(2)$.
In terms of the Woronowicz basis generators, the generators of $U_q(su(2))$
in the Drinfel'd-Jimbo basis are (see Appendix A):
\be
 J_3 & = & \frac{1}{4 \ln q} \ln \left( 1-(q-q^{-1}T_3)\right) \\
 J_{\pm} & = & T_{\pm}\left(1-(q-q^{-1}T_3\right)^{1/4}q^{\mp 1/2}.
\ee
Thus the $q$-rotations and boosts can be written in this basis as
\be
 \cj_1 & = & \frac{1}{2} \left( T^i_+(1-(q-q^{-1})T_3^i)^{1/4} q^{-1/2}
	+ T^i_-(1-(q-q^{-1})T_3^i)^{1/4} q^{1/2} \right) \\
 \cj_2 & = & \frac{1}{2i} \left( T^i_+(1-(q-q^{-1})T_3^i)^{1/4} q^{-1/2}
	- T^i_-(1-(q-q^{-1})T_3^i)^{1/4} q^{1/2} \right) \\
 \cj_3 & = & \frac{1}{4 \ln q} \ln \left[(1-(q-q^{-1})T_3^1)
			 (1-(q-q^{-1})T_3^2)\right] \\
 \ck_1 & = & -\frac{i}{2} \left( (-1)^{i-1} T^i_+(1-(q-q^{-1})T_3^i)^{1/4}
		     q^{-1/2}
	+ (-1)^{i-1}T^i_-(1-(q-q^{-1})T_3^i)^{1/4} q^{1/2} \right) \\
 \ck_2 & = & -\frac{1}{2} \left( (-1)^{i-1}T^i_+(1-(q-q^{-1})T_3^i)^{1/4}
		     q^{-1/2}
	- (-1)^{i-1}T^i_-(1-(q-q^{-1})T_3^i)^{1/4} q^{1/2} \right) \\
 \ck_3 & = & -\frac{i}{4 \ln q} \ln \frac{1-(q-q^{-1})T_3^1}
			 {1-(q-q^{-1})T_3^2},
\ee
where $i=1,2$ is summed.

In either formulation these generators
satisfy the following commutation relations:
\be
 \left[\cj_i,\cj_j\right] & = & \left\{\begin{array}{lll}
     \frac{1}{2}i\ep_{ijk}\left[2\cj_k\right],    &  \mbox{if  } & k\neq 3 \\
     \frac{1}{2}i\ep_{ijk}\left\{2\cj_k\right\}_q, &  \mbox{if  } & k=3
				      \end{array} \right. \\
 \left[\ck_i,\ck_j\right] & = & \left\{\begin{array}{lll}
     -\frac{1}{2}i\ep_{ijk}\left[2\cj_k\right],    &  \mbox{if  } & k\neq 3 \\
     -\frac{1}{2}i\ep_{ijk}\left\{2\cj_k\right\}_q, &  \mbox{if  } & k=3
				      \end{array} \right. \\
 \left[\cj_i,\ck_j\right] & = & \left\{\begin{array}{lll}
     \frac{1}{2}i\ep_{ijk}\left[2\ck_k\right],    &  \mbox{if  } & k\neq 3 \\
     \frac{1}{2}i\ep_{ijk}\left\{2\ck_k\right\}_q, &  \mbox{if  } & k=3
				      \end{array} \right. ,
\ee
where
\be
 \left\{2\cj_k\right\}_q & = & \left[\frac{1}{i}(i\cj_k-\ck_k)\right]_q
	  + \left[\frac{1}{i}(i\cj_k+\ck_k)\right]_q \\
 \left\{2\ck_k\right\}_q & = & -i\left(\left[\frac{1}{i}(i\cj_k-\ck_k)\right]_q
	  - \left[\frac{1}{i}(i\cj_k+\ck_k)\right]_q \right),
\ee
and $[A]_q$ represents the normal
q-integers. We see that in the limit $q \rightarrow 1$ the generators $\cj_i$
and $\ck_i$ reduce to the standard Lorentz rotations and boosts.

In this formalism the map between the $q$-SHO operators and the $q$-Lorentz
generators is formally identical to the classical case
and the $q$-structure is explicit only in the
deformed commutation relations.
This suggests that,
by construction, the generators of a $q$-spinorial representation,
$\si^{\mu\nu}$, should have a matrix representation in terms of the
rotation and boost generators equivalent to the classical case,
with the $q$-structure hidden in its algebra. This immediately implies a
deformed Dirac (Clifford) algebra for $\ga$-matrices. The form of this
deformation has not been determined, but we note in the next section that
the form of an associated deformed Minkowski metric can be obtained
in a simple manner from the $q$-boson operators.

\section{{$q$}-Minkowski Space}
The existence and form of a $q$-deformed Minkowski spacetime have been
considered in detail recently by a number of authors
\cite{schmidke91,ogiev92,kulish93,azcar94,podles94,azcar95}. Here we shall
simply consider how the deformed Minkowski space algebra of Ogievetsky et al.
can be related to the algebra of the $q$-SHO.

In previous discussions \cite{schmidke91,ogiev92} generators
for $q$-Minkowski space coordinates were
built out of $q$-spinor bilinears, their commutation relations being determined
by the $SU_q(2)$ $q$-spinor relations. In these discussions
the $q$-spinors were taken as
coordinates of the non-commutative quantum plane $\{(x,y):xy=qyx\}$,
the underlying carrier space of $SU_q(2)$ \cite{song92a}.
However, we now show that the $q$-Minkowski
generators can be constructed directly out of the $q$-boson operators.
Consider the
basis $\{A,B,C,D\}$, where the generators are defined as follows:

\be
 A & = & q^{-(N_{q1}-N_{q2})} \\
 B & = & q^{-1/2(N_{q1}-N_{q2})} a_1^{\dagger}a_2 \\
 C & = & \frac{(q-q^{-1})^2}{q} a_2^{\dagger}a_1 q^{-1/2(N_{q1}-N_{q2})} \\
 D & = & \frac{(q-q^{-1})^2}{q} a_2^{\dagger}a_1a_1^{\dagger}a_2
	     + q^{(N_{q1}-N_{q2})},
\ee
where $\{a_i^{\dagger},a_i,N_{qi}\}$, $i=1,2$ are
the generators of two commuting
$q$-SHOs. We note that this map differs from the construction in terms of
$q$-oscillators, and $SU_q(2)$ generators, considered by Kulish
\cite{kulish93}.
This operator basis generates the $q$-Minkowski algebra of Schmidke et al.
\cite{schmidke91}, in the basis of Kulish \cite{kulish93},
\be
 AC & = & q^2 CA \nonumber\\
 AB & = & q^{-2}BA \nonumber\\
 AD & = & DA \label{qm}\\
 \left[B,D\right] & = & -\frac{q-q^{-1}}{q}AB \nonumber\\
 \left[C,D\right] & = & \frac{q-q^{-1}}{q}CA \nonumber\\
 \left[B,C\right] & = & \frac{q-q^{-1}}{q}(AD-A^2) \nonumber,
\ee
where central elements are: the $q$-trace,
\be
 \mbox{tr}_q(K) & = & \mbox{Tr} I_q K = \frac{1}{q}A + qD \\
 I_q & = & \left(\begin{array}{cc}
    q^{-1}  &   0  \\
    0       &   q
	       \end{array}\right);
\ee
and the invariant Minkowski length,
\be
 L & = & CB - \frac{1}{q^2}DA \\
 0 & = & \left[L,\{A,B,C,D\}\right].
\ee
Kulish showed, via consideration of the reflection equation, that the
homomorphism,
\be
 K = \left(\begin{array}{cc}
     A & B \\
     C & D
	   \end{array}\right) & \rightarrow & M K \tilde{M}^{-1},
\ee
where $M$ is an element of the quantum group algebra $SL_q(2,{\boldmath C})$
(``dual'' to our formulation of $SU_q(2)$),
and $M^{\dagger} = \tilde{M}^{-1}$, leaves the defining relations
(Eq.s \ref{qm}) invariant. Thus this procedure does indeed give an analog of
Minkowski space in the $q$-deformed regime.
The $q$-Minkowski 4-vector, in $2 \times 2$ formalism, can be written in terms
of real Minkowski coordinates $\{X^0,X^1,X^2,X^3\}$ as
follows:
\be
 X & = & \left(\begin{array}{cc}
      A & B \\
      C & D
	       \end{array}\right) = \left(\begin{array}{cc}
	\frac{\sqrt{2}q}{2}(X^0-X^3)  &  \frac{\sqrt{2}q}{1-i}(X^1-X^2)\\
	\frac{\sqrt{2}}{(1+i)q}(X^1+X^2) & \frac{\sqrt{2}}{2q}(X^0+X^3)
					  \end{array}\right).
\ee
This then allows the invariant Minkowski length to be represented in the form
\be
 L & = & (q^2+1)^{-1} g_{ij} X^iX^j,
\ee
where the deformed Minkowski metric is given by
\be
 g_{ij} & = & \left(\begin{array}{cccc}
   0    &   q^2   &  0    &   0      \\
   1    &   0     &  0    &   0      \\
   0    &   0     &  0    &  -1      \\
   0    &   0     & -1    & q(q-q^{-1})
		    \end{array}\right),
\ee
as was obtained by Ogievetsky et al. \cite{ogiev92}. Clearly such a metric
does not reduce to the standard form in the $q \rightarrow 1$ limit.

\section{Concluding Remarks}
In this letter we have considered various maps from SHO operators to the
generators of $SU(2)$ and the Lorentz group. Consideration of equivalent maps
in the $q$-deformed regime has led to a $q$-deformation of the algebra
of Lorentz rotations and boosts, and
its corresponding chirally decomposed form.
We have also noted that $q$-deformed Minkowski 4-vectors can also be obtained
via construction from the $q$-SHO. The relationship between these algebras
suggests that the $q$-SHO is fundamentally related not only to the
compact quantised universal enveloping algebras
but also to the structure of the
$q$-deformed spacetime symmetry groups.
The fact that the $q$-SHO has also been found to
possess a quantum group (Hopf algebra) structure \cite{yan90,flor91,oh95}
implies that it could play a central role in the classification of
quantum groups.

It was noted earlier that the deformation of the rotation and boost algebra
of the Lorentz group should lead to a deformation of the $\ga$-matrix
Clifford algebra associated with the spinorial representation. This was
not considered here but we note in passing that such deformed Clifford
algebras have been obtained \cite{carow91,song92b} using $R$-matrix methods,
and consideration of the differential calculus on quantum spaces has allowed
Song \cite{song92b} to determine the form of the Dirac equation for
$q$-deformed Dirac bi-spinors.

We finally note that that the analysis here has focussed on the spacetime
symmetry applications of the various mappings. If, however, we were to
concentrate on the compact direct product group $SO(4) \cong SU(2) \otimes
SU(2)$ then the analysis presented here could easily be modified to  allow
consideration of the $q$-deformation of
internal symmetries such as chiral symmetry, i.e. $SU_{qL}(2)
\otimes SU_{qR}(2)$.

\appendix

\section{Comparison of {$SU_q(2)$} Generators in Various Bases}
Quantum group theory has developed from a number of rather different starting
points. The mathematical structure has been abstracted from fields as
seemingly diverse as quantum inverse scattering theory, solvable two
dimensional statistical models, non-commutative geometry, not to mention
knot theory. The manifestation of this structure in such diverse fields
is certainly suggestive of deeper underlying relationships. However, one
problem associated with the multitude of approaches to the subject is that
there have been a number of different formulations of the quantum group
structure. In particular the most widely studied $q$-deformed algebra,
$SU_q(2)$ (or $U_q(su(2))$), has been formulated in a number of different ways.
It is not immediately apparent whether these different formulations are in
fact equivalent. For discussion of the $q$-deformed Lorentz algebra it
is useful to have direct mappings between the generators of $SU_q(2)$ in
different formulations.

To achieve this we start with the following basis $\{q^{J_3},J_{\pm}\}$
of $SU_q(2)$, which
has the following algebra (Basis 1):
\be
 q^{J_3}J_{\pm}q^{-J_3} & = & q^{\pm 1} J^{\pm} \label{kul1}\\
 \left[J_+,J_-\right] & = & \frac{q^{2J_3}-q^{-2J_3}}{q-q^{-1}} \label{kul2},
\ee
due in a somewhat different form, to Kulish and Reshitikhin \cite{kulish81}.
It can easily be shown that the following $2 \times 2$ matrices constitute
a fundamental representation for the generators:
\be
 \rh(q^{J_3}) & = & \left(\begin{array}{cc}
   q^{1/2}      &    0      \\
   0            & q^{-1/2}
			  \end{array}\right) \\
 \rh(J_+) & = & \left(\begin{array}{cc}
   0      &    1      \\
   0      &    0
			  \end{array}\right) \\
 \rh(J_-) & = & \left(\begin{array}{cc}
   0      &    0      \\
   1      &    0
			  \end{array}\right).
\ee
The standard basis for the quantum universal enveloping algebra $U_q(su(2))$,
due to Drinfel'd and Jimbo \cite{drinfeld86,jimbo85,jimbo86a},
which we denote (Basis 2),
\be
 \left[J_3,J_{\pm}\right] & = & \pm J_{\pm} \label{jim1}\\
 \left[J_+,J_-\right] & = & \frac{q^{2J_3}-q^{-2J_3}}{q-q^{-1}}
	    = \frac{\sinh(2\et J_3)}{\sinh(\et)} \label{jim2},
\ee
where $q=e^{\et}$, is then a direct consequence of the algebra above
\cite{burdik91}. Eq. \ref{jim2} follows as a direct consequence of Eq.
\ref{kul2}. Eq. \ref{jim1} follows from Eq. \ref{kul1} by noting that
iteration of Eq. \ref{kul1} implies
\be
 \ln (q^{J_3}) J_{\pm} & = & J_{\pm} \ln (q^{\pm}q^{J_3}).
\ee
The fundamental representation for the generators in this basis is simply
given by a linear combination of Pauli matrices (e.g. \cite{rosso87}). i.e.
\be
 \rh(J_3) & = & \frac{1}{2}\left(\begin{array}{cc}
   1      &    0      \\
   0      &    -1
			  \end{array}\right) \\
 \rh(J_+) & = & \left(\begin{array}{cc}
   0      &    1      \\
   0      &    0
			  \end{array}\right) \\
 \rh(J_-) & = & \left(\begin{array}{cc}
   0                                   &    0      \\
   \frac{\sinh(2\et)}{\sinh(\et)}      &    0
			  \end{array}\right).
\ee

The basis for $SU_q(2)$ obtained by Woronowicz \cite{woron87,woron87b}
via analysis of non-commutative differential calculus has also been
widely studied. This basis (Basis 3) $\{T_3,T_{\pm}\}$ has the algebra
\newpage
\be
 q^{-1}T_+T_- - qT_-T_+ & = & T_3  \label{wor1}\\
 q^2T_3T_+ - q^{-2}T_+T_3 & = & (q+q^{-1}) T_+ \label{wor2} \\
 q^{-2}T_3T_- - q^2T_-T_3 & = & -(q+q^{-1})T_-, \label{wor3}
\ee
which at first sight may appear rather different to the bases of $SU_q(2)$
previously considered. However, a mapping between the generators does
exist, as was noted by Rosso \cite{rosso87}. We can obtain an explicit
mapping in two stages as follows. First create the new generator
basis $\{\ta_3,T_{\pm}\}$ via the mapping \cite{pillin94}
\be
 \ta_3 & = & q^{-4J_3} \\
 T_{\pm} & = & q^{\pm 1/2}J_{\pm}q^{-J_3}.
\ee
It then follows, from the algebra of Basis 1 (Eq.\ref{kul1},\ref{kul2}),
that
\be
 q^{-1}T_+T_- - qT_-T_+ & = & \frac{1-\ta_3}{q-q^{-1}} \\
 T_+\ta_3 & = & q^4\ta_3 T_+  \label{int2}\\
 T_-\ta_3 & = & q^{-4}\ta_3 T_-. \label{int3}
\ee
We can now obtain the algebra of the Woronowicz basis via the relation
\cite{ogiev91},
\be
 \ta_3 & = & 1 - (q-q^{-1}) T_3.
\ee
Eq. \ref{wor1} follows immediately, and Eq.s \ref{wor2},\ref{wor3} follow
from Eq.s \ref{int2},\ref{int3}. This basis has the fundamental representation
\be
 \rh(T_3) & = & \left(\begin{array}{cc}
   q^{-1}      &    0      \\
   0           &   -q
			  \end{array}\right) \\
 \rh(T_+) & = & \left(\begin{array}{cc}
   0      &    1      \\
   0      &    0
			  \end{array}\right) \\
 \rh(T_-) & = & \left(\begin{array}{cc}
   0      &    0      \\
   1      &    0
			  \end{array}\right).
\ee

For the cases of the Drinfel'd-Jimbo and Woronowicz bases (2 and 3),
which are
of interest in the discussion of the $q$-Lorentz group, we can now give
the explicit map between the generators. i.e.
\be
 T_3 & = &  \frac{1-q^{4J_3}}{q-q^{-1}} \\
 T_{\pm} & = & J_{\pm}q^{-J_3\pm 1/2}.
\ee
The inverse map is a little more complex and is given by
\be
 J_3 & = & \frac{1}{4 \ln q} \ln \left(1-(q-q^{-1})T_3 \right) \\
 J_{\pm} & = & T_{\pm}\left(1-(q-q^{-1})T_3\right)^{1/4}q^{\mp 1/2}.
\ee

\end{document}